\documentclass[letterpaper, 10 pt, conference]{ieeeconf}  % Comment this line out if you need a4paper

\IEEEoverridecommandlockouts                              % This command is only needed if 
\usepackage{graphicx} % for pdf, bitmapped graphics files
\usepackage{color}
\usepackage{amsmath} % assumes amsmath package installed
\usepackage{amssymb}  % assumes amsmath package installed
\usepackage{dsfont}
\usepackage{hyperref}

\title{\LARGE \bf
Zonal congestion management mixing large battery storage systems and generation curtailment
}

\author{Clementine Straub$^{1,2}$, Sorin Olaru$^{1}$, Jean Maeght$^{2}$ and Patrick Panciatici$^{2}$% <-this % stops a space
\thanks{*This work has been supported by the RTE-CentraleSup{\'e}lec Chair "The Digital Transformations of Electrical Networks"(https://rtechair.fr/)}
\thanks{$^{1}$ Laboratory of Signals and Systems (L2S), CentraleSup\'elec,
CNRS, Universit\'e Paris-Saclay, France. \textit{\{sorin.olaru,clementine.straub\}@centralesupelec.fr}}%
\thanks{$^{2}$ French transmission system operator, \textit{Reseau de Transport d'Electricite} (RTE), 9 rue de la porte de Buc, F-78000.}%
}

\begin{document}

\maketitle
\thispagestyle{empty}
\pagestyle{empty}

%%%%%%%%%%%%%%%%%%%%%%%%%%%%%%%%%%%%%%%%%%%%%%%%%%%%%%%%%%%%%%%%%%%%%%%%%%%%%%%%
\begin{abstract}

The French transmission system operator (RTE) needs to face a significant congestion increase in specific zones of the electrical network due to high integration of renewable energies. Network reconfiguration and renewable energy curtailment are currently employed to manage congestion and guarantee the system security and stability. In sensitive zones, however, stronger levers need to be developed. Large battery storage systems are receiving an increasing interest for their potential in congestion management. In this paper, a model for local congestion management mixing batteries and renewable generation curtailment is developed. Subsequently, an energy management approach relying on the principles of Model Predictive Control is presented. Results of simulations on RTE data sets are presented for the analysis of the degrees of freedom and sensitive parameters of the design.
\end{abstract}

%%%%%%%%%%%%%%%%%%%%%%%%%%%%%%%%%%%%%%%%%%%%%%%%%%%%%%%%%%%%%%%%%%%%%%%%%%%%%%%%
\section{INTRODUCTION}

\subsection{Motivations}

The renewable energy sector is fast-growing world-wide and the latest reports confirm such an increase in renewable energy generation on the French territory \cite{BP_RTE}. Consequently, an increase of transits in fragile zones of the electrical transportation network is foreseen. As demand is expected to stagnate, the French Transmission System Operator policy is to investigate new exploitation methods of the existing electrical installations and favour their optimal operation in the renewed context instead of developing new installations. Congestion management is a sensitive aspect in the current operation that will become critical in the future. Two means have been identified as possible technological approaches to deal with this problem: renewable generation curtailments and operation of large battery storage systems. Generation curtailment is already applied with simple strategies. For example, policies for curtailing all or half of the generation are currently implemented. These strategies are not optimal and finding the right amount to curtail, as well as the concomitant use of storage, opens the way to energy and economic savings. From a control theoretic point of view, curtailment and storage are control actions presenting delays. Their impact cannot be neglected, as the estimated delay value for curtailment is 45 seconds. In this time-delay control context \cite{olaru2008predictive}, the exploitation of the installation is facing a considerable complexity. Overloads can be allowed on some electrical lines, but are strictly regulated. The same action can clear an overload on a line and worsen flows on neighbour lines. The human operators in charge of the operation and supervision need to be assisted by an automatic policy selection for the optimal use of the new levers which, while offering new solutions to congestion problems also increase the complexity of the decision-making. An online optimization-based strategy is needed to determine what the most appropriate action is, taking into consideration the time-delay levers system, the permissible overload remaining duration and more generally network constraints.

Curtailment and storage have already been studied in several papers. For instance, \cite{gu2014fast} and \cite{burke2011factors} are investigating wind curtailment as a consequence of congested transportation networks. Taking into account that the prediction mechanism is the natural counterpart for delays \cite{normey2008dead}, the work reported in \cite{li2013mpc} presents a model predictive control (MPC) approach using battery energy systems to mitigate wind power intermittencies. The authors of \cite{parisio2014model} propose a model and optimization of microgrids bringing into account energy storage, disturbance on renewable energy generation, and curtailment schedule. However the latter articles are focused on balancing generation and load, and are not taking into account power grid limitations. The approach in \cite{biegel2012model} designs a controller for ensuring balance between power consumption and generation as well as taking into account grid capacities and thereby reducing congestion problems. The controller, though, does not take into account long delays that can occur when using curtailment and more importantly, the model is designed for tree networks. The French electrical transportation network under consideration is a mesh network. The aim of the present paper is to develop a new congestion management method designed for mesh networks that combines in a generic receding-horizon optimization problem the different levers (storage and generation curtailment) while considering their time-delay characteristics. 

It should be noted that closed-loop systems appear in the power systems literature \cite{tang2017real} \cite{bolognani2015distributed} for applications that operate on short timescales or information lacking context. The presented method needs a feedback mechanism of the system for both reasons. The controller can only have access to local measurements and needs to act fast (an action taken every two seconds). The method relies on the receding-horizon principle to reach a feedback formulation and benefits from its intrinsic prediction mechanism to deal with the time-delay in the control actions.

The main contribution of the problem is to reformulate the congestion problem in a dynamical state-space description including the operational constraints, suitable for a Model Predictive Control formulation \cite{camacho2013model} \cite{maciejowski2002predictive} \cite{mayne2000constrained}. It will be shown that the structure of the constraints is imposing a time-varying  receding-horizon optimization. The feasibility in enforced by means of a constraint softening approach and the performances are studied in a series of simulations scenarios for a representative test zone.

\subsection*{Notations}
We first define some notations.
\begin{itemize}
    \item  $\mathcal{Z}^{Nodes}$ is the set of nodes in the zone considered ; $n^{N}$ is its cardinal.
    \item  $\mathcal{Z}^{Batt} \subset \mathcal{Z}^{Nodes}$ is the set of nodes where a battery is installed ; $n^{B}$ is its cardinal.
    \item  $\mathcal{Z}^{Curt} \subset \mathcal{Z}^{Nodes}$ is the set of nodes where the generation can be curtailed ; $n^{C}$ its cardinal.
    \item  $\mathcal{Z}^{Lines} \subset \{(i,j)\in \{1..N\}^2\}$ is the set of lines in the zone considered, $n^{L}$ its cardinal.
    \item $F_{ij}$ refers to the power flow on line $ij$
    \item $P_n^{curt}$ is the amount of curtailment at node $n \in \mathcal{Z}^{Curt}$
    \item $P_k^{batt}$ is the battery power injection at node $k \in \mathcal{Z}^{Batt}$ 
    \item $E_k^{batt}$ refers to the battery energy at node $k \in \mathcal{Z}^{Batt}$ 
\end{itemize}

\section{Modeling}

\subsection{Definition of a Zone}

The French transmission electrical network is a 6000 nodes network with voltage levels from 63 to 400 kV \cite{josz2016ac}. The model considered in the present work deals with a small zone with less than 30 nodes (Fig. \ref{savignacZone}). Variations in flows are represented  by a linearization, based on report coefficient called PTDF (Power Transfer Distribution Factor).  A short description is given above and a complete definition can be found in \cite{vsovsicfeatures}. PTDF are commonly used in power systems network modeling and analysis. They are relatively easy to compute and often very useful in congestion modeling \cite{liu2002effectiveness},\cite{bart2005network}.

\subsection{Principle}

The system is modelled with the following linear representation  of the power flow equations by the mean of PTDF:

\begin{equation}
F_{ij}=F_{ij}^0+\sum_{n \in \mathcal{Z}^{Nodes}} PTDF(ij,n) \cdot \Delta P_n^{inj}
\label{PTDF}
\end{equation}
$F_{ij}^0$ is the measured power flow on line $ij$ before any system operation and $F_{ij}$ the power flow on line $ij$ after an operation. System operations considered in this article (generation curtailment and battery operation) are both injection modifications. $\Delta P_n^{inj}$ represents the injection modification at bus $n$. Equation (\ref{PTDF}) allows us to determine the new power flow resulting from an operation in a set of buses.

\begin{figure}[h]
    \centering
   \includegraphics[scale=0.8]{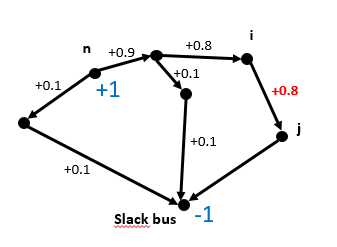}
   \caption{\label{PTDFfig} Example of PTDF on a small network}
\end{figure}

For the sake of illustration of the PTDF notion the 8 nodes benchmark network described in Fig. \ref{PTDFfig} will be used.

In order to define PTDF, a slack bus needs to be introduced in the formalism. The slack node can be any of the 8 nodes, but can not be changed along the real-time operation and in particular in the optimization-based decision making. A PTDF shows the linearized impact on a line of a transfer of power between a bus of the network and the slack bus. PTDF$(ij,n)$ is the oriented flow between $i$ and $j$ when the node $n$ produced 1 MW and the slack bus consumes 1 MW. In the example of Fig. \ref{PTDFfig}, $PTDF(ij,n)=0.8$ (in red). A definition of a PTDF is:
\begin{equation*}
    PTDF(ij,n)=
    \frac{\Delta F_{ij}}{\Delta P_n}
\end{equation*}
with $\Delta F_{ij}$ the variation in power through branch $ij$ due to the transaction from bus $n$ to the slack bus and $\Delta P_n$ the power of transaction from bus $n$ to the slack bus.

Generation curtailment and battery actions at a specific bus are modelled as a transfer of power between this bus and the slack bus. PTDF enable us to determine the fraction of power transfers flowing over each line of the zone. They are computed on the complete network with the slack bus chosen far from the zone. In a real case, $F_0$ are available via real-time measurements, while in a simulator, they are computed using the model of the complete network. PTDF and $F_0$ contain implicitly the interactions between the zone considered and the rest of the network. It may be noted that the flows representation with PTDF is equivalent to the DC modeling \cite{zimmerman2010matpower}.

\subsection{System dynamics}

The system dynamics can be represented with the help of PTDF to determine the influence of levers on flows on each line. The system outputs are curtailment and battery orders. The system presents delays between the decision-making and the order realization. The delay is $\tau^{curt} = 45s$ for generation curtailment and $\tau^{batt} = 1s$ for batteries action. The state is composed of power flows, battery charge, battery power injections and curtailments. Orders given at each time step are modifications of curtailment and power injections in the battery, respectively $\Delta P^{curt}$ and $\Delta P^{batt}$. Power injections are not constant on all time steps, they evolve according to $\Delta P^{evol}$ representing the difference between two consecutive time steps injections (changes in production or load not related to orders).

The system dynamics contains \eqref{PTDF} with

\[\Delta P^{inj} = \Delta P^{curt} + \Delta P^{batt} + \Delta P^{evol} \] 

The dynamical model has to include  the evolution of battery charges $E$, battery power injections $P^{batt}$ and curtailment $P^{curt}$, described  at each discrete-time step $k$ using a sampling of $\Delta t$ on the continuous-time scale:

\begin{equation*}
    E_b^{k+1} = E_b^k+\Delta t \cdot P_b^{batt,k+1}, \forall b \in \mathcal{Z}^{Batt}
\end{equation*}
\begin{equation*}
    P_b^{batt,k+1} = P_b^{batt,k}+\Delta P_b^{batt}, \forall b \in \mathcal{Z}^{Batt}
\end{equation*}
\begin{equation*}
    P_n^{curt,k+1} = P_n^{curt,k}+ \Delta P_n^{curt}, \forall n \in \mathcal{Z}^{Curt}
\end{equation*}

The aggregated state dynamics can be written in the form:
\begin{equation}
x^{k+1}=Ax^k+B_{curt}.u_{curt}^{k-d_{curt}}+B_{batt}.u_{batt}^{k-d_{batt}} +  B_{w}w^k
\label{eqEtat}
\end{equation} 
%\vspace{0.1cm}
with $x^k$ a vector containing power flows for each line of the zone, batteries energy, amounts of curtailed generation, and power injections in batteries:

\begin{center}
$x^k=
\begin{pmatrix}
(F_{ij}^k) \\
(E_b^k) \\
(P^{curt,k}_p) \\
(P^{batt,k}_b)
\end{pmatrix}$
$\in \mathbb{R}^n$
\end{center}
and $n={n^L+n^B+n^C+n^B}$. The discrete-time delay is related to the sampling time:
$$d_{curt}=\left\lceil \frac{\tau_{curt}}{\Delta t} \right\rceil-1;\; d_{batt}=\left\lceil \frac{\tau_{batt}}{\Delta t} \right\rceil-1$$%\vspace{0.5cm}

$u_{curt}^k$ is a vector containing orders on curtailment modifications and  $u_{batt}^k$ is a vector containing orders on modifications of power injections in batteries. {
$w^k$ is a vector representing disturbances in power injections. $w^k$ can be measured at time $k$, but $w^{k+t}$ for $t\in \{1..N\}$ is unknown. The assumption made in section \ref{simul} for predictions purposes is $w^{k+t} = w^k$ for $t \in \{0..N-1\}$} with $N$ the length of the prediction window.
$$u_{curt}^k=(\Delta P_n^{curt,k} )\in \mathbb{R}^{n^C},$$ 
$$u_{batt}^k=(\Delta P_n^{batt,k} )\in \mathbb{R}^{n^B},$$
$$ w^k=(\Delta P_n^{evol,k} )\in \mathbb{R}^{n^N} $$

The matrices $B_{curt}\in \mathbb{R}^{n\times n^C}$, $B_{batt}\in \mathbb{R}^{n\times n^B}$ and $B_{w}\in \mathbb{R}^{n\times n^N}$ are defined as follows. 

$$
B_{curt} = 
\begin{pmatrix}
M_{curt} \\
0_{n^B \times n^C} \\
\mathds{1}_{n^C\times n^C} \\
0_{n^B \times n^C}
\end{pmatrix},
B_{Batt} = 
\begin{pmatrix}
M_{batt} \\
\Delta t*\mathds{1}_{n^B \times n^B} \\
0_{n^C \times n^B} \\
\mathds{1}_{n^B \times n^B}
\end{pmatrix}$$
%{\color{blue}
$$B_{w} = 
\begin{pmatrix}
M_{w} \\
0_{(n-n^L) \times n^B}
\end{pmatrix}$$%}
%\vspace{0.5cm}

$M_{curt}$, $M_{batt}$ and $M_{w}$ contains the PTDF, such that the $k^{th}$ line in these matrices corresponds to the PTDF of the $k^{th}$ line of $(F_{ij})$ at nodes where generation can be curtailed, at nodes where a battery is installed or at nodes where the injections may vary.

All states are considered to be measured along the operation and can be initialized accordingly. Delay on the controls $u_{curt}$ and $u_{batt}$ are known on time windows $[0,d^{curt}]$ and $[0,d^{batt}]$. The system dynamics are thus well-defined in terms of forward trajectories.

\subsection{Constraints}

\label{constraints}
\begin{figure}[h]
   \includegraphics[scale=0.45]{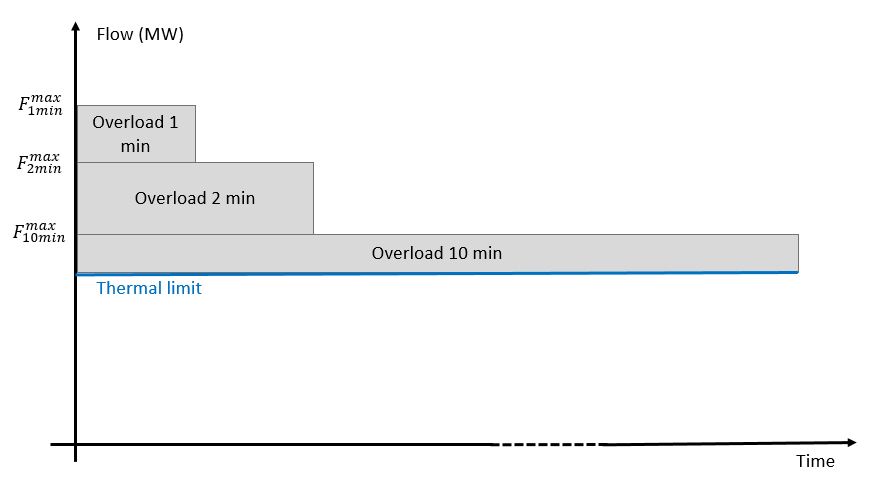}
   \caption{\label{étiquette} Example of overload capabilities}
\end{figure}

The electrical lines temperatures should stay within predefined limit values \cite{mementoRTE}. They depend on power flow (by means of the Joule effect) and weather conditions and are defined to avoid an excessive heating and dilatation. Flows must satisfy these thermal limitations. Congestion appears whenever the limit of a line is exceeded. Overloads can be tolerated for a limited time. Fig. \ref{étiquette} illustrates an example of permitted overloads, where a 15 MW overload can be tolerated for a 1 minute period.

However these overloads are allowed only when an incident occurs on the network (failure of a power plant or an electric line...) and thus should be considered as \emph{event triggered}. In a normal situation, power flows must stay within the thermal limit. As power flows are constantly evolving, a margin is taken with this limit to avoid crossing it. Assuming that control can be applied every $\Delta t = 2s$, the margin must represent the maximum of the variations during a $\Delta t$ period in normal situation.

When imposing the constraints to the system, it is necessary to distinguish between these two cases : normal situation and incident situation. Model constraints on power flow are:
\begin{equation}
|F_{ij}^{k+t}|\leq L_{ij}^{k+t}, t\in \mathbb{N}
\label{constraintsForm}
\end{equation}
$L_{ij}^{k+l}$ is defined on a finite time-horizon either by a constant function (equal to the thermal limit on line $ij$, minus the value of the margin) or a stairway profile as in Fig. \ref{étiquette}, depending on the type of situation.

If a congestion on a line cannot be eliminated in the allotted time, the line is automatically disconnected. This can result in a cascading lines opening, which is one of the main danger in power system operation. Note however, that this should not be considered as an instability from the dynamical point of view but rather as a drastic loss of performances.

The system \eqref{eqEtat} is subject to power flow limits described in \eqref{constraintsForm}. The shape of power flows limits (constant or stairway profile) depends on the network state. It is fixed before the resolution of the current optimal control problem and can only change between two resolutions of the optimal control problem. System constraints also contain bounds on batteries capacities ($E^{min}, E^{max}$ and $P^{batt,min}, P^{batt,max}$) and bounds on generation curtailment ($P^{curt,max}$). The state and input constraints can be written in the following form:
\[H_x x^k + H_u^C u_{curt}^k + H_u^B u_{batt}^k \leq H^k_0 \]
 
with

$H_x = 
\begin{pmatrix}
\mathds{1}_{n^L \times n^C} & 0 & 0 & 0 \\
0 & -\mathds{1}_{n^B \times n^B} & 0 & 0 \\
0 & \mathds{1}_{n^B \times n^B} & 0 & 0 \\
0 & 0 & \mathds{1}_{n^C \times n^C} & 0 \\
0 & 0 & 0 & -\mathds{1}_{n^B \times n^B}\\
0 & 0 & 0 & \mathds{1}_{n^B \times n^B} \\
\end{pmatrix}$,

\vspace{0.5cm}
$H_u^C = 
\begin{pmatrix}
0_{n^L \times n^C} \\
0_{n^B \times n^C} \\
0_{n^B \times n^C} \\
\mathds{1}_{n^C\times n^C} \\
0_{n^B \times n^C} \\
0_{n^B \times n^C}
\end{pmatrix}$,
$
H_u^B = 
\begin{pmatrix}
0_{n^L \times n^B} \\
\Delta t*\mathds{1}_{n^B\times n^B} \\
\Delta t*\mathds{1}_{n^B\times n^B} \\
0_{n^C \times n^B} \\
-\mathds{1}_{n^B\times n^B} \\
\mathds{1}_{n^B\times n^B}
\end{pmatrix}$,

\vspace{0.5cm}
$
H_0^k = 
\begin{pmatrix}
(L_{ij}^k)_{ij \in  \mathcal{Z}^{Lines}}\\
(E^{min})_{n \in \mathcal{Z}^{Batt}} \\
(E^{max})_{n \in \mathcal{Z}^{Batt}} \\
(P^{curt,max}_n)_{n \in \mathcal{Z}^{Curt}} \\
(P^{batt,min}_n)_{n \in \mathcal{Z}^{Batt}} \\
(P^{batt,max}_n)_{n \in \mathcal{Z}^{Batt}}
\end{pmatrix}$
\vspace{0.5cm}

The index $k$ in the matrix $H_0^k$ refers to the limits on the lines that can vary with time (constant function or stairway profile).

\subsection{Control strategy}

The system is able to receive orders, such as a desired charge level of batteries (if a battery charges to address the congestion problem, it will have to discharge at some point to be able to charge again when a new congestion problem appears). Orders can also be given regarding the generation curtailment. However, repetitive curtailment requests should be penalized. We define two objectives for the controller:
\begin{itemize}
    \item \textit{a stage cost}: $J_1(x) = \sum\limits_{t=1}^{N} ||x^{k+t}-x^{k+t}_{ref}||^2_{Q_{1}}$. The stage cost represents the desired battery charge and the desired curtailment level.
    \item \textit{a control cost}: $J_2(u) = \sum\limits_{t=0}^{N-1} ||u^{k+t}-u^{k+t}_{ref}||^2_{Q_{2}}$ representing the cost of each control. The choice of weightings allows penalizing curtailment more than battery usage. 
\end{itemize}

$||x||^2_{Q_{1}} := x^T Q_1 x$, with $Q_1\succeq 0$.

\section{Model predictive control for congestion management}

\subsection{MPC for time-delay systems}

The foregoing model offers the key elements for a prediction-based control strategy allowing the evaluation of the impact of the decisions beyond the dead-time. By considering a finite receding-horizon optimal control problem, the  following formulation can be obtained for the decision making at time step $k$:

\begin{equation}
\begin{aligned}
\label{MPCformulation}
& \underset{x,u}{\text{min}}
& & J_1(x)+J_2(u) \\
& \text{s.t.}
& & x^{k+t+1} = Ax^{k+t}+B_{curt} \cdot u_{curt}^{k+t-d^{curt}} \\
&&& \quad\quad\quad +B_{batt} \cdot u_{batt}^{k+t-d^{batt}} \\
&&& \quad\quad\quad +B_{w} \cdot w^{k+t}, \; t \in \{0,N-1\}, \\
&&& H_x x^{k+t} + H_u u^{k+t} \leq H^{k+t}_0, \; t \in \{0,N\}, \\
\end{aligned}
\end{equation}

Within the predictive control formulation, the length of the prediction window represents the main tuning parameter from both feasibility/stability and performance point of view  \cite{lombardi2012predictive, laraba2017linear}. In order to enhance the prediction and ensure its capability to cope with the time-delay, the length of the horizon must be greater than the maximal time-delay  (in this case $\tau^{curt} = 45s$ corresponding to $d^{curt}=22$): $$N\geq \max(d_{curt},d_{batt}).$$ The state estimation problem for $x^k$ is avoided as long as this information is available via measurements. However, the uncertainty has to be considered with respect to the power flow along the prediction horizon.

The matrix $H^{k+t}_0$ contains the time-dependant constraint on the overloaded lines flow (\ref{constraintsForm}). This constraint also depends on the initial flow on lines. The change in constraints from one time step to another can be handle with a reformulation of the problem as a multiparametric program. This reformulation will also deal with the problem of orders previously sent. These orders do not appear in the state $x^{k}$ and have been considered as parameters for the optimization problem. This choice is mainly related to the length of the delay and the impact on the structure of the constraints in the finite-time optimal control problem \eqref{MPCformulation}.

\subsection{MPC as a parameterized optimization problem}

The previous problem is reformulated as a parameterized optimization problem \cite{bemporad2002explicit,grancharova2012explicit}. The parameters are including the past control inputs, the current state and the measured disturbance. The constraints taken into account are time varying but their structure remains linear along one prediction window. Thus, the finite-dimensional optimization can be written as:
{
\begin{equation}
\label{paramMPC}
\begin{aligned}
& \underset{U_{batt},U_{curt}}{\text{min}}
& & J(x_k, w_k, U^P_{batt}, U^P_{curt}, U_{batt}, U_{curt}) \\
& \text{s.t.}
& & g(x_k, w_k, U^P_{batt}, U^P_{curt}, U_{batt}, U_{curt}) \leq 0\\
\end{aligned}
\end{equation}
with a quadratic cost $J(.)$, a linear set of constraints $g(.)$ and $U^P_{batt}, U^P_{curt}$ the vectors collecting the orders sent previously to the input of the system.
$$U^P_{batt} = [u_{batt}^{k-d_{batt}},u_{batt}^{k-d_{batt}+1},..., u_{batt}^{k-1}],$$
$$U^P_{curt} = [u_{curt}^{k-d_{curt}},u_{curt}^{k-d_{curt}+1},..., u_{curt}^{k-1}],$$
$$U_{batt} = [u_{batt}^{k},u_{batt}^{k+1},...,u_{batt}^{k+N-1}],$$
$$U_{curt} = [u_{curt}^{k},u_{curt}^{k+1},...,u_{curt}^{k+N-1}]$$}

The optimal solution $U_{batt}^{*}(x_k, w_k, U^P_{batt}, U^P_{curt})$ and $U_{curt}^{*}(x_k, w_k, U^P_{batt}, U^P_{curt})$ can be computed online efficiently based on Quadratic Programming (QP) solvers as the optimization problem is convex. For the simulations presented in section~\ref{simul}. the solver Fico~Xpress \cite{Xpress} was used .

\subsection{Feasibility}

The recursive feasibility of the optimization problem cannot be guaranteed using the classical arguments \cite{mayne2000constrained} which employ the tail of the optimal sequence at the previous time instant. Indeed, the structure of the constraints in \eqref{constraintsForm} is time-varying and the switch in the structure is activated by the exogenous signal $w_k$ and thus a priori unknown in the prediction scheme. For instance, gaps between power flows and lines capabilities can be too important to be addressed, or sudden increases in power flow cannot be reduced immediately because of 
input delays. The feasibility can be studied with two methods: the elimination of time-dependence or the introduction of slack variables. The elimination of time-dependence in the formulation (\ref{paramMPC}) is computationally prohibitive and over-conservative from the performance viewpoint as it relies on the enumeration of all the possible structural constraints and is therefore not employed here. The second method is to introduce slack variables $\epsilon \in \mathcal{Z}^{Lines}$ on the constraint (\ref{constraintsForm}):
\[|F_{ij}^{k+t}|\leq L_{ij}^{k+t} + \epsilon_{ij}^{k+t}, t \in \{1..N\}\]
\[\epsilon_{ij}^{k+l} \geq 0 \]
and 
$\epsilon$ is included in the cost function by adding the penalizing term:
\[J_s = \sum_{t=1}^{N} \sum_{ij \in \mathcal{Z}^{Lines}} (\epsilon_{ij}^{k+t})^2\]

The problem is then feasible for each time $t$ at the price of overloads whenever $\epsilon>0$. It is obvious to see that it exists $\epsilon \geq 0$ to render $U_{batt} = 0$ and $U_{curt} = 0$ as an admissible solution.

\section{Simulations}
\label{simul}
\begin{figure}[h]
   \centering
   \includegraphics[scale=0.45]{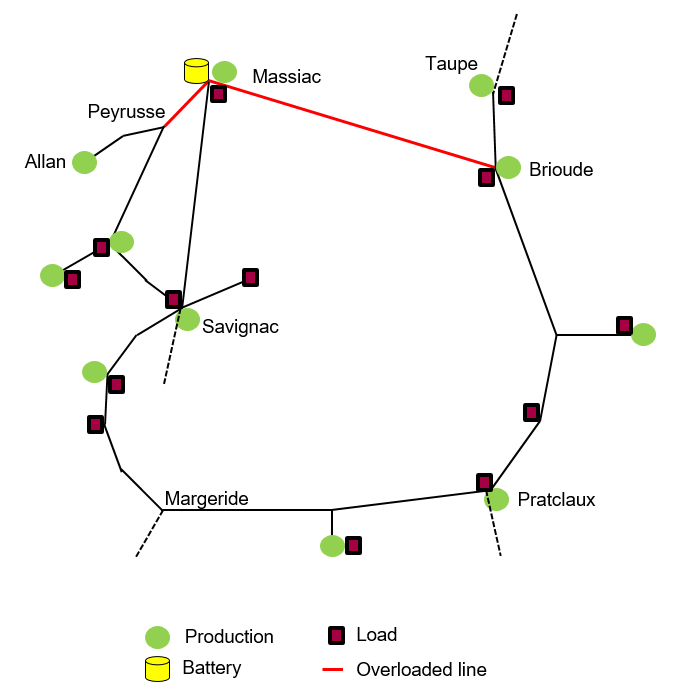}
   \caption{\label{savignacZone} Map of Savignac zone}
\end{figure}

The behaviour of the MPC scheme has been simulated on a zone around Savignac in the Auvergne region (center of France), see the map of the zone in Fig. \ref{savignacZone}. The dimensions  of the problem (\ref{MPCformulation}) are $n^{N} = 21$, $n^{B} = 1$, $n^{C} = 11$, $n^L = 22$. It follows that dim$(x^k) = 35$, dim$(u_{curt}^k) = 22$ and  dim$(u_{batt}^k) = 1$.

A major growth in renewable generation is expected in this area and congestions are anticipated. The introduction of a battery is currently studied. For our simulations, we situate the battery in Massiac and consider a capacity of $30MW - 30MWh$. The section shows the behaviour of the MPC controller facing three cases:
\begin{itemize}
    \item Constant loads and generation with one overloaded line
    \item Constant loads and generation with two overloaded lines
    \item Volatile injections
\end{itemize}

In the three cases, overloads appear due to incidents on the network. This means as described in \ref{constraints} that flows can cross thermal limits, but should stay below the stairway profile describing permitted overloads. We show that in the three cases, flows respect the permitted overloads, thanks to the MPC controller action. 

Data used has been extracted from RTE network analysis tool Convergence\footnote{See \cite{josz2016ac} for Convergence software short description}. Limitations and injections have been modified in the data, so that congestion can be observed. The modifications are representative of the increase in power flow due to the growth of renewable energy. Data for the first two cases has been simplified to better illustrates the behaviour of the controller. Generation and loads are kept constant once the incident has occurred, meaning that flows without any action of the controller are constant. A drop in a flow on a line is thus the result of the controller action only and not just the result of a variation in generation or loads.

%Power flows measurement every two seconds are not available. Data present measurements every hour. These data are nevertheless considered as if they were every two seconds. The example presented below shows the results of the MPC scheme.

\subsection*{Constant loads and generation with one overloaded line}

\begin{figure}[h!]
    \centering
   \includegraphics[scale=0.35]{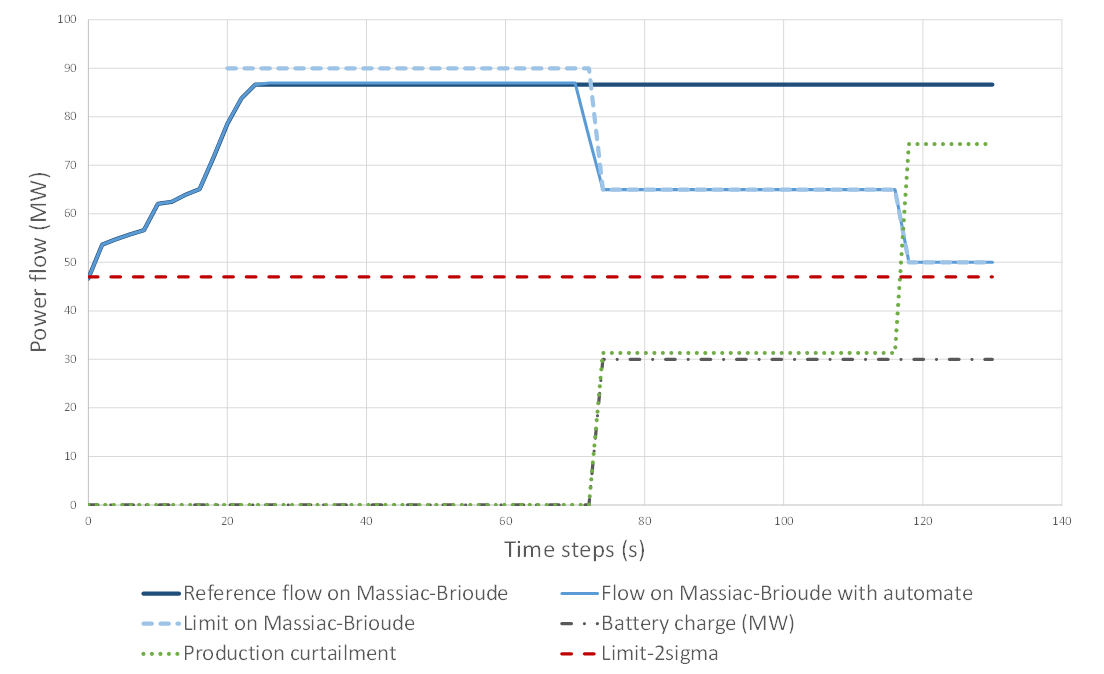}
   \caption{\label{1ligne} Power flows on Massiac-Brioude with constant generation and load, one overload}
\end{figure}

\begin{figure}[h!]
    \centering
   \includegraphics[scale=0.35]{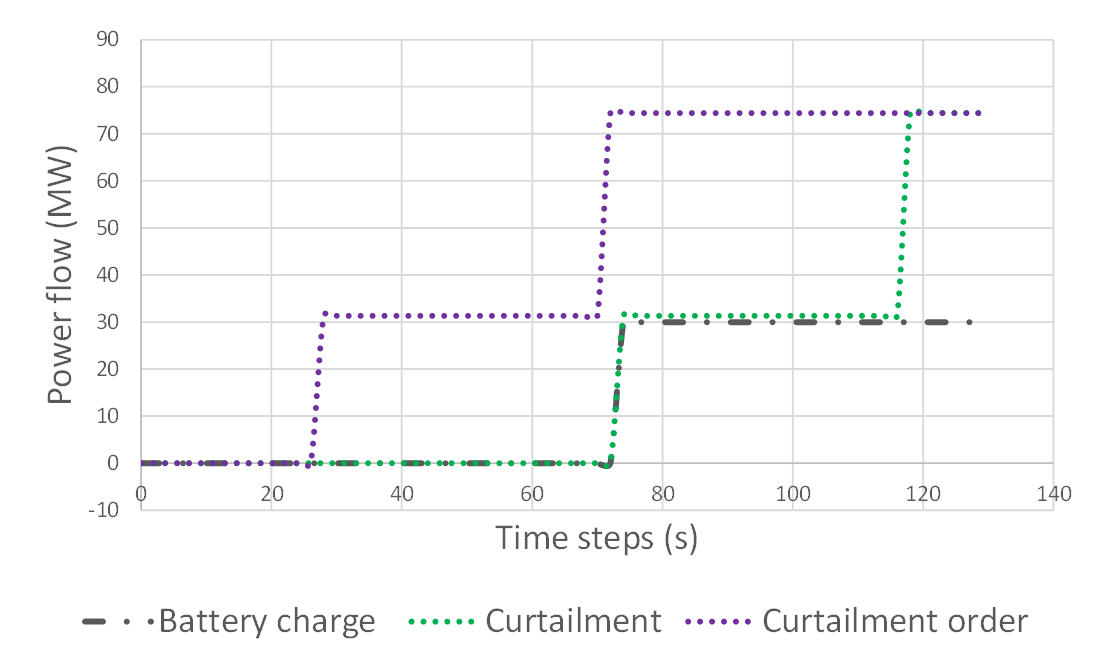}
   \caption{\label{controle1ligne} Control with constant generation and load, one overload}
\end{figure}

Fig. \ref{1ligne} and Fig. \ref{controle1ligne} illustrate the behaviour of the controller when one line is congested: the line between Massiac and Brioude. Fig. \ref{1ligne} shows the flow on the congested line with and without the controller action, as well as the actions taken. Fig. \ref{controle1ligne} illustrates the control: when the actions are taken and when they are carried out.

An incident occurred at time 0, meaning that the flow should stay under the stairway profile describing permitting overloads. The flow increases progressively until $t=24s$ and then is stabilizing around 88 MW. The reference flow (without any action of the controller) is settled at 88 MW after 24s. This flow is violating capacity constraints. The graph shows the difference between the reference flow and the flow on Massiac-Brioude when the controller is acting and when it is not. The controlled flow respects these constraints: the battery charges at its maximal power (30 MW) and generation is curtailed to bring the flow down. The controller can curtail generation in 11 places: it chose Massiac because of the values of the PTDF (see TABLE \ref{tablePTDF}). The biggest PTDF is in Massiac ; for the same amount of curtailment, the decrease in Massiac-Brioude flow will be more important if the curtailment is in Massiac in comparison to other places. The battery, located in Massiac, is used at its maximum power as it is cheaper than curtailment. The curtailment order in Massiac is sent at $t=27s$, $46s$ before the order realization at $t=72s$ : the optimization problem shows that the battery action is not sufficient for guaranteeing the feasibiblity of the congestion problem and a complementary lever is used.
The other lines are considered to be far from their limit and are not represented in the graph for readability. 

The strategy of the controller here is simple to analyze: the battery is used in priority as the cheapest lever and generation with the largest impact (PTDF) is then curtailed. Optimal curtailment comes as an automatic decision with the resolution of the optimization problem.

\subsection*{Constant loads and generation with two overloaded lines}

\begin{figure}[h!]
   \centering
   \includegraphics[scale=0.35]{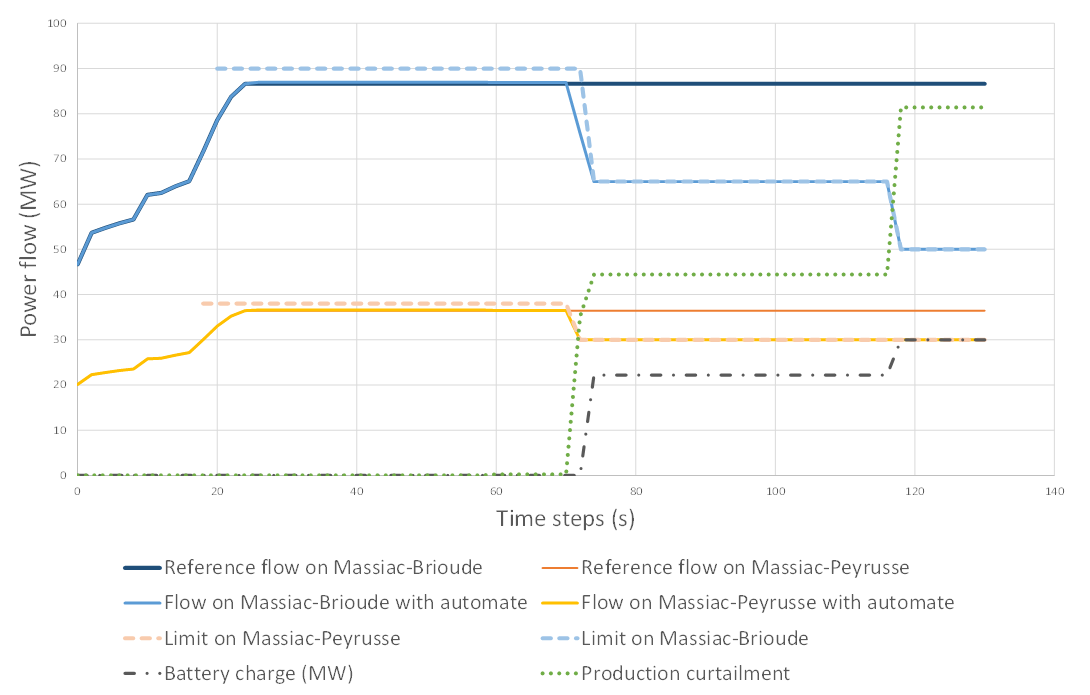}
   \caption{\label{2lignes} Power flows on Massiac-Brioude with constant generation and load, two overloads}
\end{figure}

\begin{figure}[h!]
    \centering
   \includegraphics[scale=0.35]{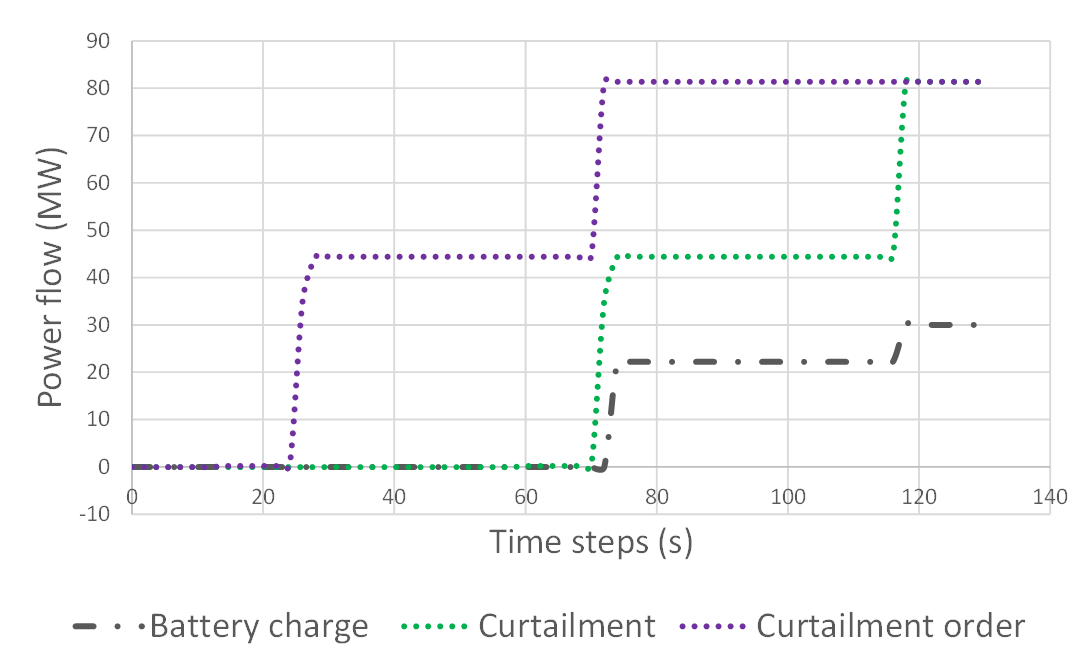}
   \caption{\label{controle2lignes} Control with constant generation and load, two overloads}
\end{figure}

Fig. \ref{2lignes} and Fig. \ref{controle2lignes} show a situation where several overloads appear at the same time. Massiac-Peyrusse line is now congested, in addition to Massiac-Brioude. As before, we see on the graph that the controller acts to reduce the flows on the overloaded lines in order to respect the capacity constraints. Both flows on Massiac-Brioude and Massiac-Peyrusse decreases with the controller action. The graph shows that the battery charges 23 MW (instead of 30 MW in the previous case), curtails more generation and not in the same place: the curtailment is here in Allan.  It can be explain by the values of the PTDF (TABLE \ref{tablePTDF}). The PTDF for Massiac on the lines Massiac-Brioude and Massiac-Peyrusse have an opposite sign: an action in Massiac will have an opposite effect on Massiac-Brioude and Massiac-Peyrusse lines. Curtailment and battery charge in Massiac only is not an admissible solution as it will increase the flow on Massiac-Peyrusse. PTDF in Allan have the same sign: curtailment in Allan decreases both the flows on Massiac-Brioude and Massiac-Peyrusse. It is the optimal solution of the optimization problem. We observe that the battery is also charging: this action decreases the flow on Massiac-Brioude and increases the flow on Massiac-Peyrusse. The drop needed on Massiac-Peyrusse is less important than the one needed on Massiac-Brioude: curtailment in Allan decreases the flow on Massiac-Peyrusse more than necessary and the battery action, by decreasing the flow on Massiac-Brioude to stay under the limit, increases the flow on Massiac-Peyrusse, while respecting the limitation.

\begin{table}
\begin{center}
    \begin{tabular}{|l|c|c|c|}
      \hline
       & Massiac & Allan & Brioude \\
      \hline
      Massiac-Brioude & 0.36 & 0.3 & 0.32\\
      \hline
      Massiac-Peyrusse & -0.27 & 0.48 & 0.14\\
      \hline
    \end{tabular}
      \caption{\label{tablePTDF}PTDF for Massiac and Allan on the overloaded lines }
\end{center}
\end{table}

The strategy of the controller is more complex to analyze than in the previous case. The controller arbitrates between the two levers to find the optimal solution. The arbitration between the different levers depends on several factors: the PTDF values, the initial flows values, the activated capacity constraints and the levers costs in the objective function. The notion of delays has not been mentioned here because of the constant loads and generation but is presented in the next paragraph.

\subsection*{Volatile injections}

\begin{figure}[h!]
   \centering
   \includegraphics[scale=0.35]{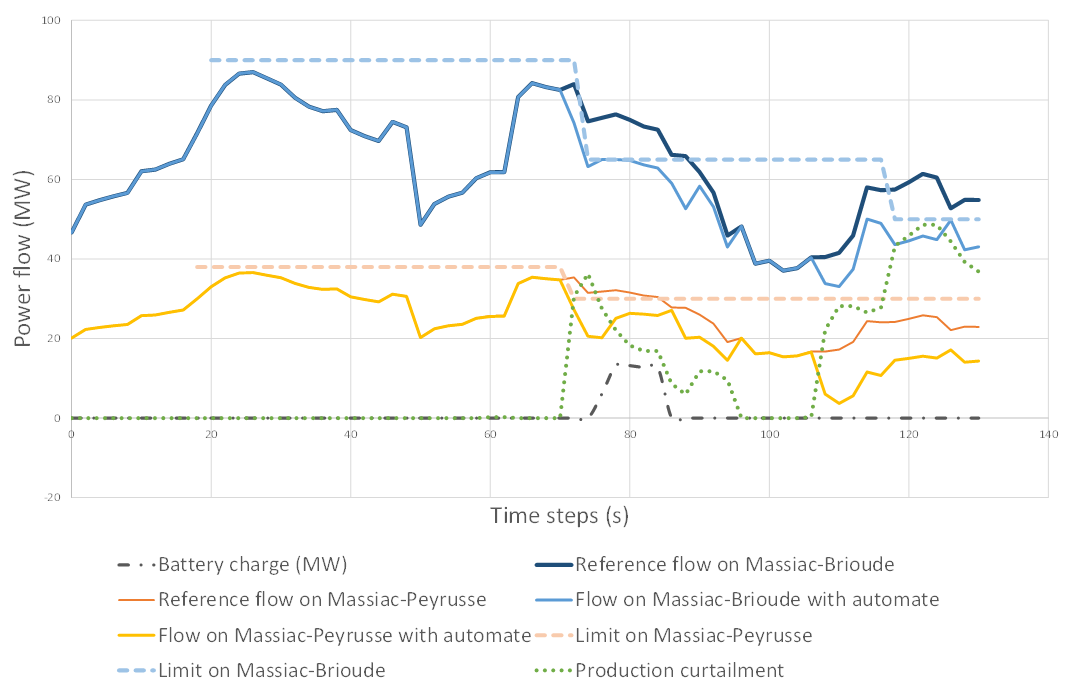}
   \caption{\label{2lignesVar} Power flows on Massiac-Brioude with volatile generation and load}
\end{figure}

\begin{figure}[h!]
    \centering
   \includegraphics[scale=0.35]{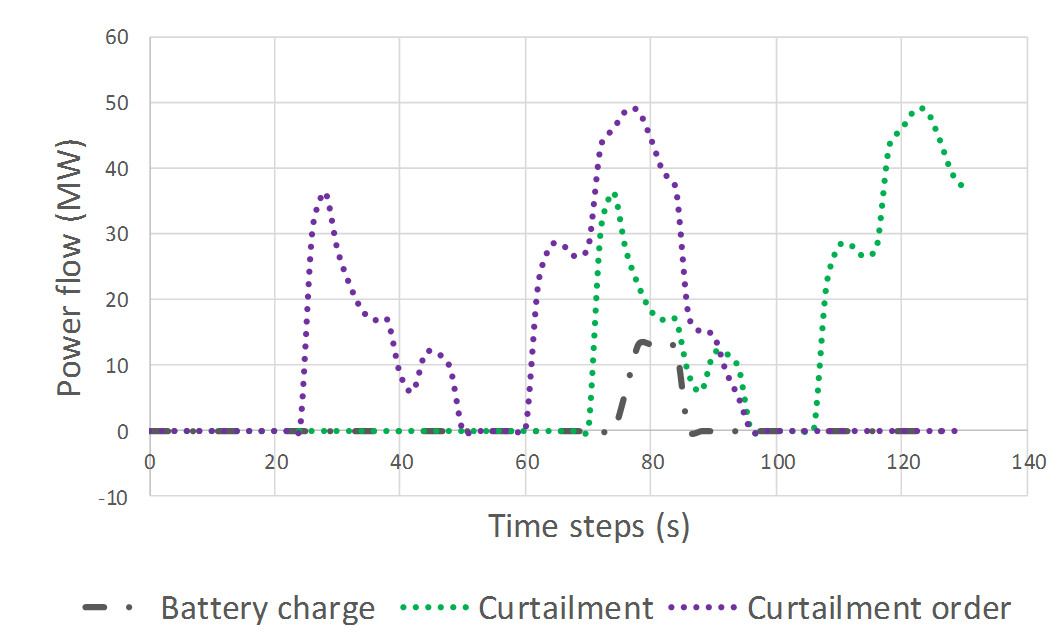}
   \caption{\label{controle2lignesVar} Control with volatile generation and load}
\end{figure}

Overloads in Fig. \ref{2lignesVar} and Fig. \ref{controle2lignesVar} case appear on the same lines: Massiac-Brioude and Massiac-Peyrusse. The other lines, far from their limits, are not represented. The two graphs illustrate the behaviour of the system with volatile injections.

 Fig. \ref{2lignesVar}  shows that flows on the two overloaded lines respect their capacity constraints as a result of the controller action. Levers used are the same as in Fig. \ref{2lignes}: curtailment in Allan and battery charge, the constraints being the same. The delays influence can be noted here. The system can not foresee variations in flows due to variations in load and generation. Its prediction horizon is constant and the assumption on the power injections ($w_{k+t} = w_k$) does not hold. As curtailment is decided 45 seconds before its real impact on the system, the amount of generation curtailed is not always optimal because of the variations between the time of the decision making and order realization. The battery is used to compensate the errors in curtailment as its action is not delayed. The controller begins to send curtailment orders at $t=26s$. The first order is effective at $t=72s$. We can note that no generation is curtailed between $t=96s$ and $t=106s$. It results from the decrease in the reference flows between $t=50s$ and $t=60s$ which are below limitations. The battery charges between $t=78s$ and $t=84s$ to compensate the insufficient curtailment (sent between $t=32s$ and $t=38s$).\\

This section illustrates the behaviour of the MPC controller and shows that levers used depend strongly on the situation of overloaded lines and on the overload on each line. On the positive side, the use of levers is accurate within considered scenario. On the negative side, the performance is strongly related to the assumptions made regarding power injections. Two options being available to mitigate this problem : either increase the prediction capability on disturbances or accept a decrease in performance by considering several scenarios and optimizing the worst case.

\section{Future work}

\subsection{Uncertainties}

Several uncertainties are present in the system and were not taking into account. It is essential to consider a system robust to these uncertainties. If a flow exceeds the capacity line, it result in an activation of the line protection : the line is opened automatically. It can disturb the electrical network operation in a dramatic way. The uncertainties are due to:

\begin{itemize}
    \item \textit{The absence of load and generation forecast.} The controller prediction horizon is constant. It is relatively short and it is not possible in practice to have precise forecasts on a short notice. For the controller, flows can vary only with control (curtailment and batteries) over the horizon, which means flows variations due to other factors (demand, wind, sun...) are not considered for the predicted trajectory.
    
    \item \textit{The approximation of the DC modelization.} It results in prediction errors.
    
    \item \textit{Unknown topology outside the observed zone}. PTDF depend on the network topology. The topology can be changed upon network operation. If PTDF are not updated, it will introduce an error in the controller flows representation.
    
    \item \textit{Imperfect measures}. Additive disturbances in the prediction model.
\end{itemize}

\subsection{Network reconfiguration}

Network reconfiguration is another possible lever. Network reconfiguration for congestion management has been studied. \cite{granelli2006optimal} presents a methodology to find the optimal topological configuration of a power transmission system using genetic algorithms. \cite{arya2000line} presents an interactive line switching algorithm for overload alleviation. Discrete variables are introduced for topological changes. Adding network reconfiguration as one of the possible level in the MPC controller will lead to an hybrid formulation. This formulation will necessitate the resolution of a Mixt Integer Programming problem every two seconds. Moreover, PTDF used in the formulation depend on the topology. Modification of the topology will lead to PTDF variations. The model will have to deal with this problem.

\subsection{Discrete generation curtailment}

As of now, the curtailment of a renewable farm is \emph{all or nothing}. It is impossible to curtail only a part of the generation and the technical interface with Distribution System Operators is difficult to change. It could be necessary to introduce discrete generation curtailment in the model. Wind farms will be curtailed taking into account not only their position but also their size.

\section{Conclusion}

We have presented a MPC controller for congestion management in electrical transmission network using high power batteries along with renewable generation curtailment. The MPC strategy proved to be well suited for this problem as it takes into account the  permitted overloads and their different duration, as well as the delays and orders on the reference charge batteries. The MPC scheme allows us to exploit the batteries speed of action and to combine them with the slower acting lever: generation curtailment. Simulation results show the behaviour of the controller in representative scenarios.

\addtolength{\textheight}{-12cm}   % This command serves to balance the column lengths
\vspace{2cm}

\bibliographystyle{plain}
\bibliography{mybib}

%%%%%%%%%%%%%%%%%%%%%%%%%%%%%%%%%%%%%%%%%%%%%%%%%%%%%%%%%%%%%%%%%%%%%%%%%%%%%%%%

%%%%%%%%%%%%%%%%%%%%%%%%%%%%%%%%%%%%%%%%%%%%%%%%%%%%%%%%%%%%%%%%%%%%%%%%%%%%%%%%

\end{document}